\documentclass[10pt]{article}
\usepackage{amsmath}
\usepackage{amssymb}
\usepackage{amsthm}
\usepackage{mathrsfs}
\usepackage[compress,noadjust]{cite}
\numberwithin{equation}{section}

\begin{document}
\title{Boundedness of singular integral operators with variable kernels on weighted weak Hardy spaces}
\author{Hua Wang \footnote{E-mail address: wanghua@pku.edu.cn.}\\
\footnotesize{Department of Mathematics, Zhejiang University, Hangzhou 310027, China}}
\date{}
\maketitle

\begin{abstract}
Let $T_\Omega$ be the singular integral operator with variable kernel $\Omega(x,z)$. In this paper, by using the atomic decomposition theory of weighted weak Hardy spaces, we will obtain the boundedness properties of $T_\Omega$ on these spaces, under some Dini type conditions imposed on the variable kernel $\Omega(x,z)$. 
\end{abstract}

\section{Introduction}

Let $S^{n-1}$ be the unit sphere in $\mathbb R^n$ ($n\ge2$) equipped with the normalized Lebesgue measure $d\sigma$. A function $\Omega(x,z)$ defined on $\mathbb R^n\times\mathbb R^n$ is said to belong to $L^{\infty}(\mathbb R^n)\times L^r(S^{n-1})$, $r\ge1$, if it satisfies the following conditions:

(1) for all $\lambda>0$ and $x,z\in\mathbb R^n$, $\Omega(x,\lambda z)=\Omega(x,z)$;

(2) for any $x\in\mathbb R^n$, $\int_{S^{n-1}}\Omega(x,z')\,d\sigma(z')=0$;

(3) $\|\Omega\|_{L^\infty(\mathbb R^n)\times L^r(S^{n-1})}:=\sup_{x\in\mathbb R^n}\left(\int_{S^{n-1}}|\Omega(x,z')|^r\,d\sigma(z')\right)^{1/r}<\infty$,

\noindent where $z'=z/{|z|}$ for any $z\in\mathbb R^n\backslash\{0\}$. Set $K(x,z)=\frac{\Omega(x,z')}{|z|^n}$. In this paper, we consider the singular integral operator with variable kernel which is defined by
\begin{equation}
T_{\Omega}f(x)=\mbox{\upshape{P.V.}}\int_{\mathbb R^n}K(x,x-y)f(y)\,dy.
\end{equation}

In \cite{cal1} and \cite{cal2}, Calder\'on and Zygmund investigated the $L^p$ boundedness of singular integral operators with variable kernels. They found that these operators $T_\Omega$ are closely related to the problem about second order elliptic partial differential equations with variable coefficients. We will denote the conjugate exponent of $p>1$ by $p'=p/{(p-1)}$. In \cite{cal3}, Calder\'on and Zygmund proved the following theorem.

\newtheorem*{thma}{Theorem A}

\begin{thma}[\cite{cal3}]
Let $1<p,r<\infty$ satisfy

$(i)$ $\frac1r<\frac1{p'}+\frac1{p'(n-1)}$ if $1<p\le2;$ or

$(ii)$ $\frac1r<\frac1{p'}+\frac{1}{p(n-1)}$ if $2\le p<\infty$.

\noindent Suppose that $\Omega(x,z)\in L^\infty(\mathbb R^n)\times L^r(S^{n-1})$. Then there exists a constant $C>0$ independent of $f$ such that
\begin{equation*}
\|T_{\Omega}(f)\|_{L^p}\le C\|f\|_{L^p}.
\end{equation*}
In particular, $T_{\Omega}$ is bounded on $L^p(\mathbb R^n)$ for all $p\ge r'$.
\end{thma}

In 1971, Muckenhoupt and Wheeden \cite{muckenhoupt2} studied the weighted norm inequalities for $T_{\Omega}$ with power weights. In 2008, Lee et al. \cite{lee1} also considered the weighted boundedness of $T_{\Omega}$ with more general weights, and showed that if the kernel $K(x,y)$ satisfies the $L^r$-H\"ormander condition with respect to $x$ and $y$ variables respectively, then $T_{\Omega}$ is bounded on $L^p_w(\mathbb R^n)$. More precisely, they proved

\newtheorem*{thmb}{Theorem B}
\begin{thmb}[\cite{lee1}]
Let $1<r<\infty$. Suppose that $\Omega(x,z)\in L^\infty(\mathbb R^n)\times L^r(S^{n-1})$ such that the following two inequalities
\begin{equation}
\sup_{x\in\mathbb R^n\atop 0<|y|<R}\sum_{k=1}^\infty\big(2^kR\big)^{n/{r'}}
\bigg(\int_{2^kR\le|z|<2^{k+1}R}\big|K(x,z-y)-K(x,z)\big|^r\,dz\bigg)^{1/r}<\infty
\end{equation}
and
\begin{equation}
\sup_{x,y\in\mathbb R^n\atop 0<|x-y|<R}\sum_{k=1}^\infty\big(2^kR\big)^{n/{r'}}
\bigg(\int_{2^kR\le|z|<2^{k+1}R}\big|K(x,z)-K(y,z)\big|^r\,dz\bigg)^{1/r}<\infty
\end{equation}
hold for all $R>0$. If $r'\le p<\infty$ and $w\in A_{p/{r'}}$, then $T_{\Omega}$ is bounded on $L^p_w(\mathbb R^n)$.
\end{thmb}

It should be pointed out that the above $L^r$-H\"ormander conditions on the variable kernels was also considered by  Rubio de Francia, Ruiz and Torrea in \cite{rubio}.

In \cite{ding6,ding7}, Ding et al. introduced some definitions about the variable kernel $\Omega(x,z)$ when they studied the $H^1$--$L^1$ boundedness of Marcinkiewicz integral. Replacing the condition (3) mentioned above, they strengthened it to the condition

$(3')$\, $\sup_{x\in\mathbb R^n\atop\rho\ge0}\left(\int_{S^{n-1}}|\Omega(x+\rho z',z')|^r\,d\sigma(z')\right)^{1/r}<\infty$.

For $r\ge1$, a function $\Omega(x,z)$ is said to satisfy the $L^r$-Dini condition if the conditions $(1)$, $(2)$, $(3')$ hold and
\begin{equation}
\int_0^1\frac{\omega_r(\delta)}{\delta}\,d\delta<\infty,
\end{equation}
where $\omega_r(\delta)$ is the integral modulus of continuity of order $r$ of $\Omega$ defined by
\begin{equation*}
\omega_r(\delta):=\sup_{x\in\mathbb R^n\atop\rho\ge0}\bigg(\int_{S^{n-1}}\sup_{y'\in S^{n-1}\atop |y'-z'|\le\delta}\big|\Omega(x+\rho z',y')-\Omega(x+\rho z',z')\big|^rd\sigma(z')\bigg)^{1/r}.
\end{equation*}

In order to obtain the $H^p_w$--$L^p_w$ boundedness of $T_{\Omega}$, Lee et al. \cite{lee1} generalized the $L^r$-Dini condition by replacing (1.4) to the following stronger condition (see also \cite{lin})
\begin{equation}
\int_0^1\frac{\omega_r(\delta)}{\delta^{1+\alpha}}\,d\delta<\infty,\quad 0\le\alpha\le1.
\end{equation}
If $\Omega$ satisfies (1.5) for some $r\ge1$ and $0\le\alpha\le1$, we say that it satisfies the $L^{r,\alpha}$-Dini condition. For the special case $\alpha=0$, it reduces to the $L^r$-Dini condition. For $0\le\beta<\alpha\le1$, if $\Omega$ satisfies the $L^{r,\alpha}$-Dini condition, then it also satisfies the $L^{r,\beta}$-Dini condition. We thus denote by $\mbox{Din}^r_\alpha(S^{n-1})$ the class of all functions which satisfy the $L^{r,\beta}$-Dini condition for all $0<\beta<\alpha$.

\newtheorem*{thmc}{Theorem C}
\begin{thmc}[\cite{lee1}]
Let $0<\alpha\le1$ and $n/{(n+\alpha)}<p<1$. Suppose $\Omega\in\mbox{Din}^r_\alpha(S^{n-1})$ such that $(1.2)$ and $(1.3)$ hold for a certain large number $r>1$. If $w^{r'}\in A_{(p+\frac{p\alpha}{n}-\frac1r)r'}$, then there exists a constant $C>0$ independent of $f$ such that
\begin{equation*}
\|T_\Omega(f)\|_{L^p_w}\le C\|f\|_{H^p_w}.
\end{equation*}
\end{thmc}

It is easy to check that
\begin{equation*}
\mbox{Din}^r_\alpha(S^{n-1})\subset\mbox{Din}^s_\alpha(S^{n-1}), \quad \mbox{if}\;\; 1\le s<r<\infty.
\end{equation*}
Then for $0<\alpha\le1$, we define
\begin{equation}
\mbox{Din}^{\infty}_\alpha(S^{n-1})=\bigcap_{r\ge1}\mbox{Din}^r_\alpha(S^{n-1}).
\end{equation}

The main purpose of this article is to study the corresponding estimates of $T_\Omega$ on the weighted weak Hardy spaces $WH^p_w(\mathbb R^n)$ (see Section 2 for its definition). We now present our main result as follows.

\newtheorem{theorem}{Theorem}[section]

\begin{theorem}
Let $0<\alpha\le1$, $n/{(n+\alpha)}<p\le1$ and $w\in A_{p(1+\frac{\alpha}{n})}$. Suppose $\Omega\in\mbox{Din}^{\infty}_\alpha(S^{n-1})$ such that $(1.2)$ and $(1.3)$ hold. Then there exists a constant $C>0$ independent of $f$ such that
\begin{equation*}
\|T_\Omega(f)\|_{WL^p_w}\le C\|f\|_{WH^p_w}.
\end{equation*}
\end{theorem}

\section{Notations and preliminaries}

The definition of $A_p$ class was first used by Muckenhoupt \cite{muckenhoupt1}, Hunt, Muckenhoupt and Wheeden \cite{hunt}, and Coifman and Fefferman \cite{coifman} in the study of weighted
$L^p$ boundedness of Hardy-Littlewood maximal functions and singular integrals. Let $w$ be a nonnegative, locally integrable function defined on $\mathbb R^n$; all cubes are assumed to have their sides parallel to the coordinate axes.
We say that $w\in A_p$, $1<p<\infty$, if
\begin{equation*}
\left(\frac1{|Q|}\int_Q w(x)\,dx\right)\left(\frac1{|Q|}\int_Q w(x)^{-1/{(p-1)}}\,dx\right)^{p-1}\le C, \quad \mbox{for every cube}\; Q\subseteq \mathbb
R^n,
\end{equation*}
where $C$ is a positive constant which is independent of the choice of $Q$. For the case $p=1$, $w\in A_1$, if
\begin{equation*}
\frac1{|Q|}\int_Q w(x)\,dx\le C\cdot\underset{x\in Q}{\mbox{ess\,inf}}\,w(x), \quad \mbox{for every cube}\;Q\subseteq\mathbb R^n.
\end{equation*}
The smallest value of $C>0$ such that the above inequalities hold is called the $A_p$ characteristic constant of $w$ and denoted by $[w]_{A_p}$.
For the case $p=\infty$, $w\in A_\infty$ if it satisfies the $A_p$ condition for some $1<p<\infty$.

A weight function $w$ is said to belong to the reverse H\"older class $RH_s$ if there exist two constants $s>1$ and $C>0$ such that the following reverse H\"older inequality holds
\begin{equation*}
\left(\frac{1}{|Q|}\int_Q w(x)^s\,dx\right)^{1/s}\le C\left(\frac{1}{|Q|}\int_Q w(x)\,dx\right), \quad \mbox{for every cube}\; Q\subseteq \mathbb R^n.
\end{equation*}

It is well known that if $w\in A_p$ with $1<p<\infty$, then $w\in A_r$ for all $r>p$, and $w\in A_q$ for some $1<q<p$. We thus write $q_w\equiv\inf\{q>1:w\in A_q\}$ to denote the critical index of $w$. Moreover, if $w\in A_p$ with $1\le p<\infty$, then there exists $s>1$ such that $w\in RH_s$. It follows directly from H\"older's inequality that $w\in RH_r$ implies $w\in RH_s$ for all $1<s<r$.

Given a cube $Q$ and $\lambda>0$, $\lambda Q$ stands for the cube with the same center as $Q$ whose side length is $\lambda$ times that of $Q$. $Q=Q(x_0,r)$ denotes the cube centered at $x_0$ with side length $r$. For a weight function $w$ and a measurable set $E$, we denote the Lebesgue measure of $E$ by $|E|$ and set the weighted measure $w(E)=\int_E w(x)\,dx$.

We give the following results that will be used in the sequel.

\newtheorem{lemma}[theorem]{Lemma}
\begin{lemma}[\cite{garcia2}]
Let $w\in A_q$ with $q\ge1$. Then, for any cube $Q$, there exists an absolute constant $C>0$ such that
$$w(2Q)\le C\,w(Q).$$
In general, for any $\lambda>1$, we have
$$w(\lambda Q)\le C\cdot\lambda^{nq}w(Q),$$
where $C$ does not depend on $Q$ or $\lambda$.
\end{lemma}

Given a weight function $w$ on $\mathbb R^n$, for $0<p<\infty$, we denote by $L^p_w(\mathbb R^n)$ the weighted space of all functions $f$ satisfying
\begin{equation}
\|f\|_{L^p_w}=\left(\int_{\mathbb R^n}|f(x)|^pw(x)\,dx\right)^{1/p}<\infty.
\end{equation}
When $p=\infty$, $L^\infty_w(\mathbb R^n)$ will be taken to mean $L^\infty(\mathbb R^n)$, and
\begin{equation}
\|f\|_{L^\infty_w}=\|f\|_{L^\infty}=\underset{x\in\mathbb R^n}{\mbox{ess\,sup}}\,|f(x)|.
\end{equation}
We also denote by $WL^p_w(\mathbb R^n)$ the weighted weak $L^p$ space which is formed by all
measurable functions $f$ satisfying
\begin{equation}
\|f\|_{WL^p_w}=\sup_{\lambda>0}\lambda\cdot w\big(\big\{x\in\mathbb R^n:|f(x)|>\lambda \big\}\big)^{1/p}<\infty.
\end{equation}

Let us now turn to the weighted weak Hardy spaces. The (unweighted) weak $H^p$ spaces have first appeared in the work of Fefferman, Rivi\`ere and Sagher \cite{cfefferman}, which are the intermediate spaces between two Hardy spaces through the real method of interpolation. The atomic decomposition characterization of weak $H^1$ space on $\mathbb R^n$ was given by Fefferman and Soria in \cite{rfefferman}. Later, Liu \cite{liu1} established the weak $H^p$ spaces on homogeneous groups for the whole range $0<p\le1$. The corresponding results related to $\mathbb R^n$ can be found in \cite{lu}. For the boundedness properties of some operators on weak Hardy spaces, we refer the readers to \cite{ding1,ding2,ding3,ding4,ding5,liu2,tao}. In 2000, Quek and Yang \cite{quek} introduced the weighted weak Hardy spaces $WH^p_w(\mathbb R^n)$ and established their atomic decompositions. Moreover, by using the atomic decomposition theory of $WH^p_w(\mathbb R^n)$, Quek and Yang \cite{quek} also obtained the boundedness of Calder\'on-Zygmund type operators on these weighted spaces.

We write $\mathscr S(\mathbb R^n)$ to denote the Schwartz space of all rapidly decreasing infinitely differentiable functions and $\mathscr S'(\mathbb R^n)$ to denote the space of all tempered distributions, i.e., the topological dual of $\mathscr S(\mathbb R^n)$. Let $w\in A_\infty$, $0<p\le1$ and $N=[n(q_w/p-1)]$. Define
\begin{equation*}
\mathscr A_{N,w}=\Big\{\varphi\in\mathscr S(\mathbb R^n):\sup_{x\in\mathbb R^n}\sup_{|\alpha|\le N+1}(1+|x|)^{N+n+1}\big|D^\alpha\varphi(x)\big|\le1\Big\},
\end{equation*}
where $\alpha=(\alpha_1,\dots,\alpha_n)\in(\mathbb N\cup\{0\})^n$, $|\alpha|=\alpha_1+\dots+\alpha_n$, and
\begin{equation*}
D^\alpha\varphi=\frac{\partial^{|\alpha|}\varphi}{\partial x^{\alpha_1}_1\cdots\partial x^{\alpha_n}_n}.
\end{equation*}
For any given $f\in\mathscr S'(\mathbb R^n)$, the grand maximal function of $f$ is defined by
\begin{equation*}
G_w f(x)=\sup_{\varphi\in\mathscr A_{N,w}}\sup_{|y-x|<t}\big|(\varphi_t*f)(y)\big|.
\end{equation*}
Then we can define the weighted weak Hardy space $WH^p_w(\mathbb R^n)$ by $WH^p_w(\mathbb R^n)=\big\{f\in\mathscr S'(\mathbb R^n):G_w f\in WL^p_w(\mathbb R^n)\big\}$. Moreover, we set $\|f\|_{WH^p_w}=\|G_w f\|_{WL^p_w}$.

\begin{theorem}[\cite{quek}]
Let $0<p\le1$ and $w\in A_\infty$. For every $f\in WH^p_w(\mathbb R^n)$, there exists a sequence of bounded measurable functions $\{f_k\}_{k=-\infty}^\infty$ such that

$(i)$ $f=\sum_{k=-\infty}^\infty f_k$ in the sense of distributions.

$(ii)$ Each $f_k$ can be further decomposed into $f_k=\sum_i b^k_i$, where $\{b^k_i\}$ satisfies

\quad $(a)$ Each $b^k_i$ is supported in a cube $Q^k_i$ with $\sum_{i}w\big(Q^k_i\big)\le c2^{-kp}$, and $\sum_i\chi_{Q^k_i}(x)\le c$. Here $\chi_E$ denotes the characteristic function of the set $E$ and $c\sim\big\|f\big\|_{WH^p_w}^p;$

\quad $(b)$ $\big\|b^k_i\big\|_{L^\infty}\le C2^k$, where $C>0$ is independent of $i$ and $k\,;$

\quad $(c)$ $\int_{\mathbb R^n}b^k_i(x)x^\alpha\,dx=0$ for every multi-index $\alpha$ with $|\alpha|\le[n({q_w}/p-1)]$.

Conversely, if $f\in\mathscr S'(\mathbb R^n)$ has a decomposition satisfying $(i)$ and $(ii)$, then $f\in WH^p_w(\mathbb R^n)$. Moreover, we have $\big\|f\big\|_{WH^p_w}^p\sim c.$
\end{theorem}

Throughout this article $C$ always denotes a positive constant, which is independent of the main parameters and not necessarily the same at each occurrence.

\section{Proof of Theorem 1.1}

Following the same arguments as in the proof of Lemma 5 in \cite{kurtz}, we can also establish the following lemma on the variable kernel $\Omega(x,z)$ (See \cite{ding7} and \cite{lee1}).

\begin{lemma}
Let $r\ge1$. Suppose that $\Omega(x,z)\in L^\infty(\mathbb R^n)\times L^r(S^{n-1})$ satisfies the $L^r$-Dini condition in Section $1$. If there exists a constant $0<\gamma\le 1/2$ such that $|y|<\gamma R$, then for any $x_0\in\mathbb R^n$, we have
\begin{equation*}
\begin{split}
\bigg(&\int_{R\le|x|<2R}\big|K(x+x_0,x-y)-K(x+x_0,x)\big|^rdx\bigg)^{1/r}\\
&\le C\cdot R^{-n/{r'}}\bigg(\frac{|y|}{R}+\int_{|y|/{2R}}^{|y|/R}\frac{\omega_r(\delta)}{\delta}d\delta\bigg),
\end{split}
\end{equation*}
where the constant $C>0$ is independent of $R$ and $y$.
\end{lemma}

We are now in a position to give the proof of Theorem 1.1.

\begin{proof}[Proof of Theorem 1.1]
For any given $\lambda>0$, we may choose $k_0\in\mathbb Z$ such that $2^{k_0}\le\lambda<2^{k_0+1}$. For every $f\in WH^p_w(\mathbb R^n)$, then by Theorem 2.2, we can write
\begin{equation*}
f=\sum_{k=-\infty}^\infty f_k=\sum_{k=-\infty}^{k_0} f_k+\sum_{k=k_0+1}^\infty f_k:=F_1+F_2,
\end{equation*}
where $F_1=\sum_{k=-\infty}^{k_0} f_k=\sum_{k=-\infty}^{k_0}\sum_i b^k_i$, $F_2=\sum_{k=k_0+1}^\infty f_k=\sum_{k=k_0+1}^\infty\sum_i b^k_i$ and $\{b^k_i\}$ satisfies $(a)$--$(c)$ in Theorem 2.2. Then we have
\begin{equation*}
\begin{split}
&\lambda^p\cdot w\big(\big\{x\in\mathbb R^n:|T_{\Omega}(f)(x)|>\lambda\big\}\big)\\
\le\,&\lambda^p\cdot w\big(\big\{x\in\mathbb R^n:|T_{\Omega}(F_1)(x)|>\lambda/2\big\}\big)+\lambda^p\cdot w\big(\big\{x\in\mathbb R^n:|T_{\Omega}(F_2)(x)|>\lambda/2\big\}\big)\\
=\,&I_1+I_2.
\end{split}
\end{equation*}
First we claim that the following inequality holds:
\begin{equation}
\big\|F_1\big\|_{L^2_w}\le C\cdot\lambda^{1-p/2}\big\|f\big\|^{p/2}_{WH^p_w}.
\end{equation}
In fact, since supp\,$b^k_i\subseteq Q^k_i=Q\big(x^k_i,r^k_i\big)$ and $\big\|b^k_i\big\|_{L^\infty}\le C 2^k$ according to Theorem 2.2, then it follows directly from Minkowski's inequality that
\begin{equation*}
\begin{split}
\big\|F_1\big\|_{L^2_w}&\le\sum_{k=-\infty}^{k_0}\sum_i\big\|b^k_i\big\|_{L^2_w}\\
&\le\sum_{k=-\infty}^{k_0}\sum_i\big\|b^k_i\big\|_{L^\infty}w\big(Q^k_i\big)^{1/2}.
\end{split}
\end{equation*}
For each $k\in\mathbb Z$, by using the bounded overlapping property of the cubes $\{Q^k_i\}$ and the fact that $1-p/2>0$, we thus obtain
\begin{equation*}
\begin{split}
\big\|F_1\big\|_{L^2_w}&\le C\sum_{k=-\infty}^{k_0}2^k\Big(\sum_i w\big(Q^k_i\big)\Big)^{1/2}\\
&\le C\sum_{k=-\infty}^{k_0}2^{k(1-p/2)}\big\|f\big\|^{p/2}_{WH^p_w}\\
&\le C\sum_{k=-\infty}^{k_0}2^{(k-k_0)(1-p/2)}\cdot\lambda^{1-p/2}\big\|f\big\|^{p/2}_{WH^p_w}\\
&\le C\cdot\lambda^{1-p/2}\big\|f\big\|^{p/2}_{WH^p_w}.
\end{split}
\end{equation*}
Since $w\in A_{p(1+\frac{\alpha}{n})}$ and $p(1+\frac{\alpha}{n})\le1+\frac{\alpha}{n}<2$, then we have $w\in A_2$. In this case, we know that there exists a number $s>1$ such that $w\in RH_s$. More specifically, by using the sharp reverse H\"older's inequality for $A_2$ weights obtained recently in \cite{chung}, we find that for $w\in A_2$,
\begin{equation*}
w\in RH_s \quad \mbox{with}\;\; s=1+\frac{1}{2^{n+5}[w]_{A_2}}.
\end{equation*}
Observe that $\Omega\in\mbox{Din}^{\infty}_\alpha(S^{n-1})$, then we are able to find a positive number $r>1$ large enough such that $r>\max\big\{s',{(2n)}/{(n-\alpha)}\big\}$ and $\Omega\in\mbox{Din}^{r}_\alpha(S^{n-1})$. By the choice of $r$, we can easily check that $2/{r'}>1+\alpha/n\ge p(1+\alpha/n)$, which implies $w\in A_{2/{r'}}$. Hence, by using Theorem B, we know that $T_{\Omega}$ is bounded on $L^{2}_w(\mathbb R^n)$. This fact together with Chebyshev's inequality and (3.1) yields
\begin{align}
I_1&\le \lambda^p\cdot\frac{4}{\lambda^2}\big\|T_{\Omega}(F_1)\big\|^2_{L^2_w}\notag\\
&\le C\cdot\lambda^{p-2}\big\|F_1\big\|^2_{L^2_w}\notag\\
&\le C\big\|f\big\|^{p}_{WH^p_w}.
\end{align}
We now turn our attention to the estimate of $I_2$. Setting
\begin{equation*}
A_{k_0}=\bigcup_{k=k_0+1}^\infty\bigcup_i \widetilde{Q^k_i},
\end{equation*}
where $\widetilde{Q^k_i}=Q\big(x^k_i,\tau^{{(k-k_0)}/{(n+\alpha)}}(2\sqrt n)r^k_i\big)$ and $\tau$ is a fixed positive number such that $1<\tau<2$. Thus, we can further decompose $I_2$ as
\begin{equation*}
\begin{split}
I_2&\le\lambda^p\cdot w\big(\big\{x\in A_{k_0}:|T_{\Omega}(F_2)(x)|>\lambda/2\big\}\big)+
\lambda^p\cdot w\big(\big\{x\in (A_{k_0})^c:|T_{\Omega}(F_2)(x)|>\lambda/2\big\}\big)\\
&=I'_2+I''_2.
\end{split}
\end{equation*}
Let us first deal with the term $I'_2$. Since $w\in A_{p(1+\frac{\alpha}{n})}$, then by Lemma 2.1, we can deduce that
\begin{align}
I'_2&\le\lambda^p\sum_{k=k_0+1}^\infty\sum_i w\big(\widetilde{Q^k_i}\big)\notag\\
&\le C\cdot\lambda^p\sum_{k=k_0+1}^\infty\tau^{(k-k_0)p}\sum_i w\big(Q^k_i\big)\notag\\
&\le C\big\|f\big\|^{p}_{WH^p_w}\sum_{k=k_0+1}^\infty\Big(\frac{\tau}{2}\Big)^{(k-k_0)p}\notag\\
&\le C\big\|f\big\|^{p}_{WH^p_w}.
\end{align}
On the other hand, it follows immediately from Chebyshev's inequality that
\begin{equation*}
\begin{split}
I''_2&\le 2^p\int_{(A_{k_0})^c}\big|T_{\Omega}(F_2)(x)\big|^pw(x)\,dx\\
&\le 2^p
\sum_{k=k_0+1}^\infty\sum_i\int_{\big(\widetilde{Q^k_i}\big)^c}\big|T_{\Omega}\big(b^k_i\big)(x)\big|^pw(x)\,dx\\
&= 2^p
\sum_{k=k_0+1}^\infty\sum_i J^k_i.
\end{split}
\end{equation*}
Now denote $\tau^k_{i,\ell}=2^{\ell-1}\tau^{{(k-k_0)}/{(n+\alpha)}}\sqrt{n}r^k_i$, $\widetilde{Q^k_{i,\ell}}=Q\big(x^k_i,\tau^k_{i,\ell}\big)$ and
$$E^k_{i,\ell}=\big\{x\in\mathbb R^n:\tau^k_{i,\ell}\le|x-x^k_i|<2\tau^k_{i,\ell}\big\},\quad \ell=1,2,\ldots.$$
An application of H\"older's inequality gives us that
\begin{equation*}
\begin{split}
J^k_i&\le\sum_{\ell=1}^\infty\int_{E^k_{i,\ell}}\big|T_{\Omega}\big(b^k_i\big)(x)\big|^pw(x)\,dx\\
&\le\sum_{\ell=1}^\infty\bigg(\int_{E^k_{i,\ell}}w(x)\,dx\bigg)^{1-p}
\bigg(\int_{E^k_{i,\ell}}\big|T_{\Omega}\big(b^k_i\big)(x)\big|w(x)\,dx\bigg)^p.
\end{split}
\end{equation*}
Let $q=p(1+\frac{\alpha}{n})$ for simplicity. Then for any $n/{(n+\alpha)}<p\le1$ and $w\in A_q$ with $q>1$, we can easily see that $[n(q_w/p-1)]=0$. Hence, by the cancellation condition of $b^k_i\in L^\infty(\mathbb R^n)$, we get
\begin{equation*}
\begin{split}
\int_{E^k_{i,\ell}}\big|T_{\Omega}\big(b^k_i\big)(x)\big|w(x)\,dx
&=\int_{E^k_{i,\ell}}\left|\int_{Q^k_i}\Big[K\big(x,x-y\big)-K\big(x,x-x^k_i\big)\Big]b^k_i(y)\,dy\right|w(x)\,dx\\
&\le\int_{Q^k_i}\bigg\{\int_{E^k_{i,\ell}}\Big|K\big(x,x-y\big)-K\big(x,x-x^k_i\big)\Big|w(x)\,dx\bigg\}
\big|b^k_i(y)\big|\,dy\\
&\le\big\|b^k_i\big\|_{L^\infty}\big|Q^k_i\big|
\bigg(\int_{E^k_{i,\ell}}\Big|K\big(x,x-y\big)-K\big(x,x-x^k_i\big)\Big|w(x)\,dx\bigg).
\end{split}
\end{equation*}
When $y\in Q^k_i$ and $x\in\big(\widetilde{Q^k_i}\big)^c$, then a trivial computation shows that
\begin{equation}
\big|x-x^k_i\big|\ge\tau^{{(k-k_0)}/{(n+\alpha)}}\sqrt{n} r^k_i>\sqrt{n} r^k_i\ge 2\big|y-x^k_i\big|.
\end{equation}
We also observe that $w\in RH_s$ and $r>s'$, then $w\in RH_{r'}$. Using H\"older's inequality, the estimate (3.4) and Lemma 3.1, we can see that for any $y\in Q^k_i$, the integral of the above expression is dominated by
\begin{align}
&\bigg(\int_{E^k_{i,\ell}}\Big|K\big(x,x-y\big)-K\big(x,x-x^k_i\big)\Big|^r\,dx\bigg)^{1/r}
\bigg(\int_{E^k_{i,\ell}}w(x)^{r'}\,dx\bigg)^{1/{r'}}\notag\\
\le\, &C\cdot\frac{w\big(\widetilde{Q^k_{i,\ell+1}}\big)}{\big|\widetilde{Q^k_{i,\ell+1}}\big|^{1/r}}
\left(\int_{\tau^k_{i,\ell}\le|x|<2\tau^k_{i,\ell}}\Big|K\big(x+x^k_i,x-(y-x^k_i)\big)
-K\big(x+x^k_i,x\big)\Big|^r\,dx\right)^{1/r}\notag\\
\le\, &C\cdot\frac{w\big(\widetilde{Q^k_{i,\ell+1}}\big)}{\big|\widetilde{Q^k_{i,\ell+1}}\big|^{1/r}}
\cdot\Big(\tau^k_{i,\ell}\Big)^{-n/{r'}}
\left(\frac{|y-x^k_i|}{\tau^k_{i,\ell}}+
\int_{|y-x^k_i|/{2\tau^k_{i,\ell}}}^{|y-x^k_i|/{\tau^k_{i,\ell}}}\frac{\omega_r(\delta)}{\delta}\,d\delta\right)
\notag\\
\le\, &C\cdot\frac{w\big(\widetilde{Q^k_{i,\ell+1}}\big)}{\big|\widetilde{Q^k_{i,\ell+1}}\big|^{1/r}}
\cdot\Big(\tau^k_{i,\ell}\Big)^{-n/{r'}}
\left(\frac{|y-x^k_i|}{\tau^k_{i,\ell}}+\frac{|y-x^k_i|^\alpha}{(\tau^k_{i,\ell})^\alpha}\times
\int_{|y-x^k_i|/{2\tau^k_{i,\ell}}}^{|y-x^k_i|/{\tau^k_{i,\ell}}}
\frac{\omega_r(\delta)}{\delta^{1+\alpha}}\,d\delta\right)
\notag\\
\le\, & C\cdot\frac{w\big(\widetilde{Q^k_{i,\ell+1}}\big)}{\big|\widetilde{Q^k_{i,\ell+1}}\big|}
\cdot\left(\frac{1}{2^{\ell}\tau^{{(k-k_0)}/{(n+\alpha)}}}+
\Big[\frac{1}{2^{\ell}\tau^{{(k-k_0)}/{(n+\alpha)}}}\Big]^\alpha
\int_{0}^{1}\frac{\omega_r(\delta)}{\delta^{1+\alpha}}\,d\delta\right)
\notag\\
\le\, & C\cdot\left(1+\int_{0}^{1}\frac{\omega_r(\delta)}{\delta^{1+\alpha}}\,d\delta\right)
\cdot\frac{w\big(\widetilde{Q^k_{i,\ell+1}}\big)}{\big|\widetilde{Q^k_{i,\ell+1}}\big|}
\left(\frac{1}{2^{\ell}\tau^{{(k-k_0)}/{(n+\alpha)}}}\right)^\alpha.
\end{align}
Recall that $\big\|b^k_i\big\|_{L^\infty}\le C 2^k$. From the above estimate (3.5), it follows that
\begin{equation*}
\begin{split}
J^k_i&\le C\cdot2^{kp}\sum_{\ell=1}^\infty w\Big(\widetilde{Q^k_{i,\ell+1}}\Big)
\bigg(\frac{|Q^k_i|}{\big|\widetilde{Q^k_{i,\ell+1}}\big|}\bigg)^p
\left(\frac{1}{2^{\ell}\tau^{{(k-k_0)}/{(n+\alpha)}}}\right)^{\alpha p}.
\end{split}
\end{equation*}
In addition, for $w\in A_q$ with $q>1$, then we can take a sufficiently small number $\varepsilon>0$ such that $w\in A_{q-\varepsilon}$. Thus, by using Lemma 2.1 again, we finally obtain
\begin{equation*}
\begin{split}
J^k_i&\le C\cdot2^{kp}w\big(Q^k_i\big)\sum_{\ell=1}^\infty
\left(2^{\ell}\tau^{{(k-k_0)}/{(n+\alpha)}}\sqrt n\right)^{n(q-\varepsilon)-np}
\left(\frac{1}{2^{\ell}\tau^{{(k-k_0)}/{(n+\alpha)}}}\right)^{\alpha p}\\
&\le C\cdot2^{kp}w\big(Q^k_i\big)\sum_{\ell=1}^\infty
\left(2^{\ell}\tau^{{(k-k_0)}/{(n+\alpha)}}\right)^{-n\varepsilon}\\
&\le C\cdot2^{kp}w\big(Q^k_i\big)\left(\tau^{{(k-k_0)}/{(n+\alpha)}}\right)^{-n\varepsilon}.
\end{split}
\end{equation*}
Therefore
\begin{align}
I''_2&\le C\sum_{k=k_0+1}^\infty\sum_i2^{kp}\Big(\tau^{{(k-k_0)}/{(n+\alpha)}}\Big)^{-n\varepsilon}w\big(Q^k_i\big)\notag\\
&\le C\big\|f\big\|^{p}_{WH^p_w}\sum_{k=k_0+1}^\infty\Big(\tau^{{(k-k_0)}/{(n+\alpha)}}\Big)^{-n\varepsilon}\notag\\
&\le C\big\|f\big\|^{p}_{WH^p_w}.
\end{align}
Combining the above inequality (3.6) with (3.2) and (3.3), and then taking the supremum over all $\lambda>0$, we conclude the proof of Theorem 1.1.
\end{proof}

\section*{Acknowledgment}

The author would like to thank Professor J. Duoandikoetxea for pointing out the reference \cite{rubio}.

\end{document}